\documentclass[twoside,a4paper,english]{article}

\def\quality{final} %
\def\material{normal} %
\usepackage{aformat}

\usepackage{babel}

\usepackage{times}

\usepackage{epsfig,graphicx,color}

\usepackage{amsmath,amssymb,bbold}

\usepackage{mymath,myamsproof,maththm,mynumb}
\def\derpar#1{\partial_{#1}}

\usepackage{varioref}
\usepackage{xr}
\newcommand{\e}{\varepsilon}

\newcommand{\ip}[3]{\left< {#1}, {#2} \right>_{#3}}

\newcommand{\mean}[1]{\mbox{avg}\left({#1}\right)}
\newcommand{\grad}[1]{\nabla_{\!\!H^{#1}}}
\newcommand{\gradt}[1]{\nabla_{\!\!\tilde{H}^{#1}}}

\newcommand{\hh}{\widehat{h}}
\newcommand{\hk}{\widehat{k}}

\newcommand{\R}{\mathbb{R}}

\newcommand{\Z}{\mathbb{Z}}

\newcommand{\ud}{\,\mathrm{d}}

\makeatletter
\def\len{\mathop{\operator@font len}\nolimits}
\def\Len{\mathop{\operator@font Len}\nolimits}
\def\dil{\mathop{\operator@font dil}\nolimits}
\def\mod{\mathop{\operator@font mod}\nolimits}
\def\supess{\mathop{\operator@font supess}\nolimits}
\makeatother
\def\E{{\mathbb E}}

\begin{document}

\title{Sobolev--type metrics in the space of curves}
\author{
 A.~C.~G.~Mennucci \thanks{Scuola Normale Superiore, Pisa, Italy}
 \and
 A.~Yezzi  \thanks{Georgia Institute of Technology, Atlanta, USA}
 \and
 G.~Sundaramoorthi \thanks{Georgia Institute of Technology, Atlanta, USA}
}
\date{}
\maketitle

\begin{abstract}
  We define a manifold $M$ where objects $c\in M$ are curves, which we
  parameterize as $c:S^1\to \real^n$ ($n\ge 2$, $S^1$ is the circle).
  Given a curve $c$, we define the tangent space $T_cM$ of $M$ at $c$
  including in it all deformations $h:S^1\to\real^n$ of $c$.

  In this paper we study geometries on the manifold of curves,
  provided by Sobolev--type metrics $H^j$.  We study $H^j$ type
  metrics for the cases $j=1,2$; we prove estimates, and characterize
  the completion of the space of smooth curves.

  As a bonus, we prove that the Fr\'echet distance of curves (see
  \cite{Michor-Mumford}) coincides with the distance induced by the
  ``Finsler $L^\infinity$ metric'' defined in \S2.2 in
  \cite{YM:metrics04}.
\end{abstract}

\section{Introduction}
\newcommand{\C}{M}

Suppose that $c$ is an immersed curve $c:S^1\to \real^n$,
where $S^1\subset\real^2$ is the circle; we want to define 
a geometry on $M$, the space of all such immersions $c$.

The tangent space $T_cM$ of $M$ at $c$ contains all the \emph{deformations}
$h\in T_cM$  of the curve $c$, that are all the vector fields along $c$. Then,
an infinitesimal deformation of the
curve $c$ in ``direction'' $h$ will yield (on first order) the curve 
$c(u)+\varepsilon h(u)$.

For the sake of simplicity, we postpone details of the definitions (in
particular on the regularity of $c$ and $h$, and the topology on $M$)
to \S\ref{sec:spaces}.

We would like to define a \emph{Riemannian metric}
on the manifold $M$ of immersed curves: 
this means that, given two deformations $h,k\in T_cM$, we want to define
a scalar product $\langle h,k\rangle_c$, possibly dependent on $c$.
The Riemannian metric would then entail a \emph{distance} $d(c_0,c_1)$ between
the curves in $M$, defined as the infimum of the length $\Len(\gamma)$
of all smooth paths $\gamma:[0,1]\to M$ connecting $c_0$ to $c_1$.
We call \emph{minimal geodesic} a path providing the minimum of $\Len(\gamma)$
in the class of $\gamma$ with fixed endpoints.
\footnote{Note that this is an oversimplification of what we will actually do:
  compare definitions \ref{rem: d on B} and \ref{def: d F}}

At the same time, we would like to consider the curves as ``geometric
objects''; to this end, we will define the space of \emph{geometrical
  curves} $B\defeq M/\mbox{Diff}(S^1)$, that is the space of immersed
curves up to reparametrization. To this end, we will ask that the
metric (and all the energies) defined on $M$ be independent of the
parameterization of the curves.

$B$ and $M$ are the \emph{Shape Spaces} that are studied in this paper.

\bigskip

A number of methods have been proposed in Shape Analysis to
define distances between shapes, averages of shapes and optimal
morphings between shapes. At the same time, there has been much
previous work in Shape Optimization, for example Image Segmentation
via Active Contours, 3D Stereo Reconstruction via Deformable Surfaces;
in these later methods, many authors have defined Energy Functionals
$E(c)$ on curves (or on surfaces), whose minima represent the desired
Segmentation/Reconstruction;  and then utilized the Calculus of Variations to
derive curve evolutions to search minima of $E(c)$, often
referring to these evolutions as Gradient Flows.
The reference to these flows as \emph{gradient flows} implies a
certain Riemannian metric on the space of curves; but this fact has
been largely overlooked.  We call this metric $H^0$, and define it by
\[ \ip{h}{k}{H^0} \defeq \frac{1}{L} \int_{S^1} \langle h(s) , k(s)
\rangle \ud s \] where $h,k\in T_cM$, $L$ is the length of $c$, $\ud
s\defeq |\dot c(\theta)|\ud \theta$ is integration by arc-parameter,
and $\langle h(s) , k(s) \rangle$ is the usual Euclidian scalar
product in $\real^n$  (that sometimes we will also write as 
$h(s) \cdot k(s)$).

For example, the well known Geometric Heat Flow is often referred as
the \emph{gradient flow for length}: we show how it is indeed the
gradient flow w.r.t.  the $H^0$ metric.
\begin{Example}
  Let $c$ be an immersed curve, and $h$ be a deformation of $c$.
  Let the differential operator $D_s\defeq \frac1{|\dot c|}\derpar \theta$ be 
  ``the derivative  with respect to arclength''.
  Let 
  \begin{equation}
    \label{eq:len}
    \len(c)\defeq \int_{S^1} |\dot c(\theta)|\ud \theta 
  \end{equation}
  be the length of the curve;   we recall that 
  \begin{equation}\label{eq:der_L_h}
    \derparr h {\len(c)} = \int_{S^1}\langle D_s h \cdot T \rangle \ud s =
    - \int_{S^1}\langle  h \cdot D_s^2 c \rangle \ud s
  \end{equation}
  where   $T=D_s c$ is the tangent to the curve,
  $D_s^2 c$ is the curvature, intended as a vector.

  Let $C=C(\theta,t)$ be an evolving family of curves trying to
  minimize $L$: we recall moreover that the resulting
  \emph{geometric heat flow} (also known as \emph{motion by mean curvature})
  \[ \derparr t C = D_s^2 C \]
  is well defined only for positive times.

  By comparing the above flow to the definition of $H^0$, we 
  realize that this flow is the gradient descent 
  (up to a conformal term $1/\len(c)$):
  \[ \derparr t C = - \frac 1  {\len(c)} \grad 0 \len(c) \]
\end{Example}

If one wishes to have a consistent view of the geometry of the space
of curves in both Shape Optimization and Shape Analysis, then one
should use the $H^0$ metric when computing distances, averages and
morphs between shapes.

Surprisingly, $H^0$ does not yield a well define metric structure,
since the associated distance is identically zero%
\footnote{This striking fact was first described
  in~\cite{Mumford:Gibbs}; it is generalized to 
  spaces of submanifolds in \cite{Michor-Mumford04}
}.

Moreover, some simple Shape Optimization tasks are ill-defined when
using the $H^0$ metric:
\begin{Remark}
  Let $g:\real^n\to\real^k$, let 
  \[\mean{g(c)} \defeq \frac{1}{L}\int_{S^1} g(c(s)) \ud s =
  \frac{1}{L}\int_{S^1} g(c(\theta))|\dot c(\theta)|\ud \theta  ~~;\]  
  (here again $L\defeq\len(c)$); then
  \begin{eqnarray}
    \derparr h {\mean{g(c)}} &=&\frac{1}{L}\int_{S^1} \nabla g(c) h
    + g(c) \langle D_s h \cdot T \rangle \ud s - \nonumber \\ & & -
    \frac{1}{L^2}\int_{S^1} g(c)\ud s
    \int_{S^1}\langle D_s h \cdot T \rangle \ud s= \nonumber\\
    &=& \frac{1}{L}\int_{S^1}
    \nabla g(c) h  + \big( g(c) -\mean{g(c)} \big) 
      \langle D_s h \cdot T \rangle \ud s
      \label{eq:der avg c}
  \end{eqnarray}

  If the curve is in the plane, that is $n=2$, then we define the
  normal  vector   $N\perp T$ by rotating $T$ counterclockwise,
  and define scalar curvature $\kappa$ so that $D^{2}_s c= \kappa N$;
  then, integrating by parts, the above becomes
  \begin{eqnarray*}
    \derparr h {\mean{g(c)}} &=&
    \frac{1}{L}\int_{S^1}
    \derparr x g(c) h -  \big(\derparr x g(c(s)) T\big)  \langle h \cdot T \rangle 
    -\\ &  &- \big( g(c) -\mean{g(c)} \big) 
    \langle  h \cdot D^{2}_s  c\rangle \ud s=\\
    &=&
    \frac{1}{L}\int_{S^1} 
    \Big(\derparr x g(c) N - \kappa  \big( g(c)  -\mean{g(c)} \big)\Big)
    \langle h \cdot N \rangle \ud s
  \end{eqnarray*}

  Suppose now that we have a Shape Optimization functional $E$
  including a term of the form $\mean{g(c)}$; 
  let $C=C(\theta,t)$ be an evolving family of curves trying to
  minimize $E$; this flow would contain a term of the form
  \[ \derparr t C  = \dots \big( g(c(s))  -\overbar g \big) \kappa N \dots\]
  unfortunately the above flow is ill defined: it is a negative-time
  heat flow on roughly half of the curve. We present two simple examples.

  \begin{itemize}
  \item If for example $g(x)=x$, then $\mean{g(c)}=\mean{c}$ is the
    \textbf{center of mass of the curve}.
    Let us fix a target point  $v\in \real^2$. Let $E(c)\defeq
    \frac 1 2 |\mean c -v|^2$ be a functional that penalizes the distance from
    the center of mass to $v$; by direct computation
    \begin{eqnarray*}
      \derparr h E   &=& (\mean{c}-v) \cdot \derparr h {\mean c } =\\
      &=&
      \frac{1}{L}\int_{S^1} 
      \big\langle (\overbar{c}-v)\cdot 
      \big(N  -    \kappa  (c-\overbar{c}) \big) \big\rangle
      \langle h \cdot N \rangle \ud s
    \end{eqnarray*}
    (where $\overbar c=\mean c$ for simplicity);
    so we conclude that the $H^0$ gradient descent flow is
    \begin{equation}
      \derparr t C = - \grad 0 E(c) = 
      \langle(v-\overbar{c})\cdot N\rangle N - \kappa N \big\langle (c -\overbar
      c) \cdot (v-\overbar{c}) \big\rangle\label{eq:desc_c_o_m}
    \end{equation}
    (up to a part tangent to the curve).

    \begin{minipage}[c]{0.26\linewidth}
      \vskip 1pt
      {\includegraphics[width=1\linewidth]{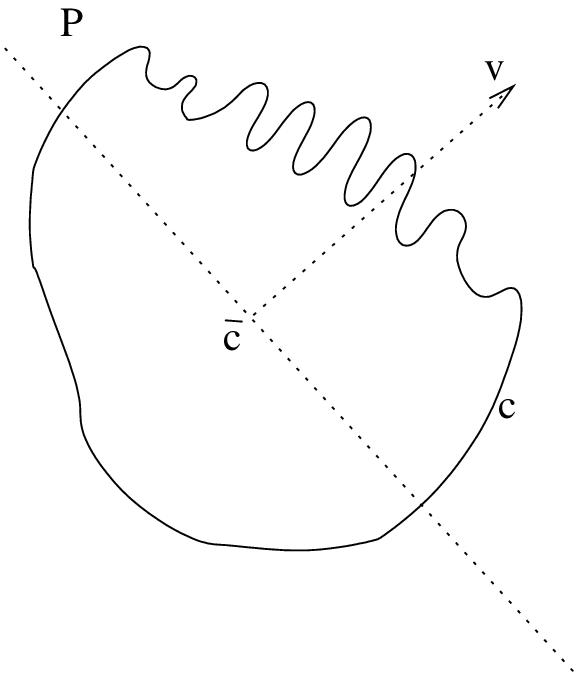}}
    \end{minipage}
    \begin{minipage}[c]{0.7\linewidth}
      Let 
      $P\defeq \{ w:\langle (w -\overbar c) \cdot (v-\overbar{c})
      \big\rangle      \ge      0\}$
      be the halfplane 
      that is the region of the plane that is ``on the $v$ side'' 
      w.r.t. $\overbar c$.
      This  gradient descent flow \eqref{eq:desc_c_o_m} does move the center of mass
      towards the point $v$: indeed there is a first term  
      $\langle(v-\overbar{c})\cdot N\rangle N$ that moves the 
      whole curve towards $v$;  and a second term that
      tries to decrease the curve length out of $P$ and
      increase the curve length
      in $P$:  and this    is ill posed.
    \end{minipage}
  
  \item 
    Similarly if $g(x)= |x-\overbar c|^2$, then $E(c)\defeq \mean{ |c
    -\overbar c|^2}$ is the \textbf{standard deviation of the curve}.
    The derivative is 
    \[\frac{2}{L}\int_{S^1} 
    \Big(  \langle(c-\overbar c) \cdot N\rangle
    - \kappa  \big( g(c)  -\mean{g(c)} \big)\Big)
    \langle h \cdot N \rangle \ud s \]

    The flow to minimize this should be 
    \[ \derparr t C = - \langle (c-\overbar c)\cdot N\rangle N + \kappa N \big( g(c)
    - \mean{g(c)} \big) \]
    and this is ill posed where the curve is inside of the 
    circle of center $\mean c$ and radius $\sqrt{\mean{g(c)}}$
  \end{itemize}
\end{Remark}

The above phenomenon is also visible in many applications, where
the Active Contour curve would ``fractalize'' in an attempt
to minimize the task energy. For this reason, a regularization term is often
added to the energy: this remedy, though, does change the energy,
and ends up solving a different problem.

\subsection{Sobolev--type Riemannian Metrics}

To overcome this limitation, in 
\cite{SYM05:sobol_activ_contour} and \cite{SYM06:track_with_sobol_activ_contour}
we proposed a family of So\-bo\-lev--type Riemannian Metrics

\begin{Definition} \label{defn:inner}
  Let $c\in \C$, $L$ be the length of $c$, and $h,k \in T_c\C$.  Let
  $\lambda > 0$. We assume $h$ and $k$ are parameterized by
  the arclength parameter of $c$. We define, for $j\ge 1$ integer,
  \begin{enumerate} 
    \addtolength{\topsep}{-1cm}
  \item $\ip{h}{k}{H^j} \defeq \ip{h}{k}{H^0} + \lambda L^{2j}
    \ip{D^j_sh}{D^j_sk}{H^0}$
  \item $\ip{h}{k}{\tilde{H}^j} \defeq \mean{h} \cdot \mean{k}+\lambda L^{2j}
    \ip{D^j_sh}{D^j_s k}{H^0}$ 
  \end{enumerate}
  where again $\mean{h} \defeq \frac{1}{L}\int_{S^1} h(s) \ud s$
  and  $D^j_s$ is the j-th  derivative with respect to arclength.
  \footnote{Note that $\ip{h}{k}{H^0}=\mean{h \cdot k}$
    so the difference in the two metrics is in using 
    $\mean{h \cdot k}$ instead of $\mean{h} \cdot \mean{k}$
  }
\end{Definition}
It is easy to verify that the above definitions are inner products.
Note that we have introduced length dependent scale factors so that the above
inner products (and corresponding norms) are independent of curve rescaling. 

Changing the metric will change the gradient and thus the gradient
descent flow; this change will alter the topology in the space of
curves, but the change of topology does not affect the energy to be
minimized, or its global minima; whereas it may regularize the flows,
and avoid that the flows be trapped in local minima; 
many examples and applications can be found in the survey \cite{ganesh:SAC}.

\subsection{Previous work}

In \cite{YM:eusipco,YM:metrics04,YM05:confor_metric}
we addressed the problem of defining a metric in the space $M$ of
parameterized immersed curves $c:S^1\to \real^n$ 
(with special attention to the case
$n=2$); we discuss the general setting of Finsler 
and Riemannian metrics, and related problematics;
we propose a set of goal properties.
we discuss some models available in the literature.
Eventually we proposed and study conformal metrics such as
\begin{equation}\label{eq:conformal}
\langle h,k\rangle_{H^0_\phi} \defeq \len(c)\int\langle h(s),k(s)\rangle ds
\end{equation}
we prove  results regarding this metric, and in particular, that the
associate distance is non degenerate. We also proved that,
supposing that only unit length curves with an upper bound on
curvature are allowed, then there exist minimal geodesics.

The same approach was proposed independently by J.~Shah in
\cite{Shah:H0}, who moreover proved that in the simplest case given by
\eqref{eq:conformal}, the only minimal geodesics are given by a curve
evolution whose velocity is proportional to the curve normal vector field.

\smallskip

The problem  has been addressed by Michor and Mumford in
\cite{Michor-Mumford}, who propose the metric
\begin{equation}
\langle h,k\rangle_{H^0_A} \defeq \int  (1+A \kappa^2(s))\langle h(s),k(s)\rangle ds
\label{eq:MM HA}
\end{equation}
where $\kappa$ is the curvature of $c$, and $A>0$ is a fixed constant;
they prove many results regarding this metric; in particular, that the
associate distance is non degenerate, and that completion of smooth
curves is contained between the space Lip of rectifiable curves, and
the space BV$^2$ of rectifiable curves whose curvature is a bounded
measure.

\medskip

More recently in 
\cite{SYM05:sobol_activ_contour} and \cite{SYM06:track_with_sobol_activ_contour}
we studied the family of Sobolev--type metrics defined in \ref{defn:inner}.

In \cite{SYM05:sobol_activ_contour} we have experimentally shown that
Sobolev flows are smooth in the space of curves, are not as dependent
on local image information as $H^0$ flows, are global motions which
deform locally after moving globally, and do not require derivatives
of the curve to be defined for region-based and edge-based energies.
In general, Sobolev gradients can be expressed in terms of the
traditional $H^0$ gradient, that is, we have the formulas
\begin{eqnarray}
  \grad{n} E  &=& \grad0 E \ast K_{\lambda,n} \\
  \label{eq:convH1t}
  \gradt{n} E &=& \grad0 E \ast \tilde{K}_{\lambda,n}
\end{eqnarray}
for suitable convolutional kernels $\tilde{K}_{\lambda,n},K_{\lambda,n}$.
We have moreover shown mathematically that the Sobolev--type gradients
regularize the flows of well known energies, by reducing the degree of
the P.D.E.
\begin{Example}
  For example, in the case of the elastic energy $ E(c) = \int_{c}
  \kappa^2 \ud s = \int_{c} |D_s^2 c|^2 \ud s $, the $H^0$ gradient is
  $\grad0 E = L D_s(2D_s^{(3)}c+3|D_s^2 c|^2D_sc)$ that includes
  fourth order derivatives; whereas the $\tilde H^1$-gradient is
  \begin{equation}
    L c_{ss} +  L  (|D_{ss}c|^2D_sc) \ast \tilde{K}_{\lambda,1}
  \end{equation}
  that is an integro-differential second order P.D.E.
\end{Example}

\smallskip

In \cite{SYM06:track_with_sobol_activ_contour} we have shown numerical
experiments on real-life cases, and shown that the regularizing
properties may be explained in the Fourier domain: indeed, if we
calculate Sobolev gradients $\grad{n} E$ of an arbitrary energy $E$ in the
frequency domain, then
\begin{equation}\label{eq:grad_n__grad_0}
  \widehat{\grad{n} E}(l) = 
  \frac{\widehat{\grad0 E}(l)}{1+\lambda(2\pi l)^{2n}} \quad
  \mbox{for $l\in \Z$}
\end{equation}
and
\begin{equation}
  \widehat{\gradt{n} E}(0) = \widehat{\grad0 E}(0) , \quad 
  \widehat{\gradt{n} E}(l) =
  \frac{\widehat{\grad0 E}(l)}{\lambda(2\pi l)^{2n}}
  \quad \mbox{for $l\in \Z \backslash \{0\}$},
\end{equation}
(see eqn.~\ref{eq:hseries} for the precise definition of Fourier coefficients).
It is clear from the previous expressions that high frequency
components of $\grad0 E(c)$ are increasingly less pronounced in the
various forms of the $H^n$ gradients.

\smallskip

\begingroup
\let\on=\operatorname
A family of metrics similar to the first example above (but for the
length dependent scale factors) is currently studied in \cite{MM:hamil_rieman}:
the Sobolev-type weak Riemannian metric on $\on{Imm}(S^1,\mathbb R^2)$
\begin{align*}
\langle h,k\rangle_{G^n_c} &=
\int_{S^1}\sum_{i=0}^n\langle D_s^ih,
   D_s^ik  \rangle ds
= \int_{S^1}\langle A_n(h),k\rangle ds
  \qquad\text{  where }
\\
A_n(h)&= A_{n,c}(h) =\sum_{i=0}^n (-1)^i D_s^{2i}(h)
\end{align*}
in that paper the geodesic equation, horizontality,
conserved momenta, lower and upper bounds
on the induced distance, and scalar curvatures are computed.
Note that this norm is equivalent to the norms in \ref{defn:inner}:
see remark~\ref{MM equiv}.
\endgroup

\bigskip

Many other approaches to Shape Analysis are present in the literature; 
for example, much earlier than the above, Younes in \cite{Younes:Comp}
had proposed a computable definition of distance of curves, modeled
on elastic curves.

\smallskip

A different approach to the study of shapes is obtained when the shape
is defined to be a \emph{curve up to reparametrization, rotation,
  translation, scaling}.

For example Mio, Srivastava et al in
\cite{Srivastava:AnalPlan,mio_sriv_04:elast_strin_model} use a
different choice of curve representation: they represent a planar
curve $c$ by a pair of angle-velocity functions $\dot c(u) =
\exp(\phi(u) + i \theta(u)) $ (identifying $\real^2=\complex$), and
then defining a metric on $(\phi,\theta)$.  They propose models of
spaces of curves where the metrics involve higher order derivatives in
\cite{Srivastava:AnalPlan}.  See
the proof of thm.~\ref{thm:H1_approx} for an example comparison
of the two approaches.

\section{Spaces of curves}
\label{sec:spaces}

As anticipated in the introduction,  we want to define 
a geometry on $M$, the space of all immersions $c:S^1\to\real^n$.

We will sometimes distinguish exactly what $M$ is, choosing between
the space $\mbox{Imm}(S^1,\real^n)$
of immersions,  the space
$\mbox{Imm}_f(S^1,\real^n)$
of \emph{free immersions},
and $\mbox{Emb}(S^1,\real^n)$
of \emph{embeddings} (see \S2.4,\S2.5 in \cite{Michor-Mumford}).

We will equip $M$ with a topology $\tau$ stronger than the $C^1$ topology:
then any such choice $M$ is an open subset of the vector space 
$C^1(S^1,\real^n)$ (that is a Banach space), so it is a manifold.

The tangent space $T_cM$ of $M$ at $c$ contains vector fields 
$h:S^1\to\real^n$ along~$c$. 

Note that we represent both curves $c\in M$ and
deformations $h\in T_cM$ as functions $S^1\to \real^n$;
this   is a special structure that is not usually present in
abstract manifolds: so we can easily define ``charts'' for $M$:
\begin{Remark}[Charts in $M$]
  \label{chart in M}
  Given a curve $c$, there is a neighbourhood $U_c$ of $0\in T_cM$
  such that for $h\in U_c$, the curve $c+h$ is still immersed; then
  this map $h\mapsto c+h$ is the simplest natural candidate to be a chart
  of $\Phi_c:U_c\to M$; indeed, if we pick another curve $\tilde c\in
  M$ and the corresponding $U_{\tilde c}$ such that $U_{\tilde c}\cap
  U_c\neq\emptyset$, then the equality $\Phi_c(h)=c+h=\tilde c+\tilde
  h=\Phi_{\tilde c}(\tilde h)$ can be solved for $h$ to obtain
  $h=(\tilde c-c)+\tilde h$.
\end{Remark}
The above  is trivial but is worth remarking for two reasons:
it stresses that the topology $\tau$ must be strong enough to mantain
immersions; and is a basis block to what we will do in the space $B_{i,f}$
defined below.

We look mainly for metrics in the space $M$ that are independent
on the parameterization of the curves $c$: to this end, 
we define these spaces of \emph{geometrical curves}
\[  B_{i} = B_{i} (S^1,\real^2) = \mbox{Imm} (S^1,\real^2)/\mbox{Diff}(S^1) \]
and 
\[  B_{i,f} = B_{i,f} (S^1,\real^2) =\mbox{Imm}_f(S^1,\real^2)/\mbox{Diff}(S^1) \]
that are the quotients
of the spaces $\mbox{Imm}_f$  and $\mbox{Imm}_f (S^1,\real^2)$
by $\mbox{Diff}(S^1)$;  alternatively we may quotient by $\mbox{Diff}^+(S^1)$
(the space of orientation preserving automorphisms of $S^1$),
and obtain spaces of \emph{geometrical oriented curves}.

\begin{Remark}[on model spaces and properties]
  We have two possible choices in mind for  the topology $\tau$ to put on $M$:
  the Fr\'echet  space of $C^\infinity$ functions; or a Hilbert space
  such as standard Sobolev space $H^j(S^1\to\real^n)$.

  Suppose we define on $M$ a Riemannian metric: we would like $B_i$ to
  have a nice geometrical structure; we would like our Riemannian Geometry to
  satisfy some useful properties.

  Unfortunately, this currently seems an antinomy.

  If $M$ is modeled on a Hilbert space $H^j$, then most of the usual
  calculus carries on; for example, the exponenential
  map would be locally a diffeomorphism; but the quotient space $M/\mbox{Diff}(S^1)$
  is not a smooth bundle, (since the tangent to the orbit contains
  $ \dot c$ and this is in $H^{j-1}$ in general!).
  
  If $M$ is modeled on the Fr\'echet space of $C^\infinity$ functions,
  then the quotient space $M/\mbox{Diff}(S^1)$ is a smooth bundle; but
  some of the usual calculus fails: the Cauchy theorem does not hold
  in general; and the  exponenential is not locally surjective.
\end{Remark}

Suppose in the following that $\tau$ is the Fr\'echet space of
$C^\infinity$ functions, for simplicity; then $B_{i,f}$ is a manifold,
the base of a principal fiber bundle while $B_{i}$ is not
(see in \S2.4.3 in \cite{Michor-Mumford}).

To define charts on this manifold, we imitate what was done for $M$:
\begin{Proposition}[Charts in $B_{i,f}$]\label{def:ortho-charts}
  Let $\Pi$ be the projection from 
  $\mbox{Imm}_f(S^1,\real^2)$ %
  to the quotient $B_{i,f}$.

  Let  $[c]\in B_{i,f}$: we pick a curve $c$ such that $\Pi(c)=[c]$.
  We represent the tangent space $T_{[c]}B_{i,f}$ as the space of all
  $k:S^1\to \real^n$ such that $k(s)$ is orthogonal to $\dot c(s)$.

  Again we can define a simple natural chart $\Phi_{[c]}$ 
  by projecting the chart $\Phi_c$ (defined in \ref{chart in M}):
  the chart is  
  \[\Phi_{[c]}(k)\defeq \Pi(c(\cdot)+k(\cdot)) \]
  that is, it moves $c(u)$ in  direction $k(u)$;
  and it is easily seen that the chart does not depend on the choice of
  $c$ such that $\Pi(c)=[c]$.
  We can solve  $\Phi_{[c]}(k)=\Phi_{[\tilde c]}(\tilde k)$
  (this is not so easy to prove: see \cite{Michor-Mumford},
  or 4.4.7 and 4.6.6 in   \cite{Hamilton:NashMoser}).
\end{Proposition}

Define a Finsler metric $F$ on $M$; this is a lower semi continuous function 
such that $F(c,\cdot)$ is a norm on $T_c M$.

If $\gamma :[0,1]\to M$ is a path connecting two curves $c_0,c_1$,
then we may define a homotopy $C:S^1\times[0,1]\to\real^n$
associated to $\gamma$ by $C(\theta,v)=\gamma(v)(\theta)$, and viceversa.

\begin{Definition}[standard distance]\label{rem: d on B}
  Given a metric $F$ in $M$, we could consequently define the 
  \emph{standard distance}  of two curves $c_0,c_1$ as the infimum of the length
  \[ \int_0^1 F(\gamma(t),\dot\gamma(t))\ud t\] in the class of all
  $\gamma$ connecting $c_0,c_1$. 
\end{Definition}
 This is not, though, the most interesting distance for applications:
 we are indeed interested in
  studying metrics and distances in the quotient space $B\defeq
  M/\mbox{Diff}(S^1)$.

We suppose that
\begin{Definition}\label{def: c w p i}
  the metric $F(c,h)$ is
  ``\textbf{curve-wise parameterization invariant}'', that is, it 
  does not depend on the parameterization of the curves $c$
\end{Definition}
then $F$ may be projected to  $B\defeq M/\mbox{Diff}(S^1)$;
we will say that $F$ is a \textbf{geometrical metric}.

\begin{minipage}[c]{0.26\linewidth}
\vskip 1pt
{\includegraphics[width=1\linewidth]{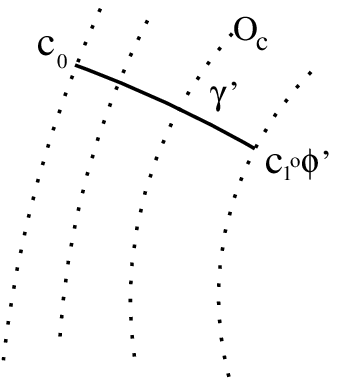}}
\end{minipage}
\begin{minipage}[c]{0.7\linewidth}
Consider two geometrical curves $[c_0],[c_1]\in B$,
and a  path $\gamma :[0,1]\to B_i$  connecting $[c_0],[c_1]$:
then we may lift it to a homotopy $C:S^1\times[0,1]\to\real^n$;
in this case, the homotopy will connect a $c_0\circ \phi_0$ to  
$c_1\circ \phi_1$, with $\phi_0,\phi_1\in \mbox{Diff}(S^1)$.
Since $F$ does not depend on the parameterization, we can factor 
out $\phi_0$ from the definition of the projected length.
\end{minipage}

To summarize, we define the
\begin{Definition}[geometric distance]\label{def: d F}
  Given $c_0,c_1$, we define the class $\mathcal A$ of homotopies $C$
  connecting the curve $c_0$ to a reparameterization 
  \footnote{If we use $\mbox{Diff}^+(S^1)$ to define $B$
    then $\phi$ must be orientation preserving as well.  }
  $c_1\circ\phi$ of
  the curve $c_1$, that is, $C(u,0)=c_0(u)$ and $C(u,1)=c_1(\phi(u))$.
  We define the \emph{geometric distance} $d_F$ of $[c_0],[c_1]$ in 
  $B\defeq M/\mbox{Diff}(S^1)$ as the
  infimum of the length
  \[\Len_F(C) \defeq \int_0^1 F\left(C(\cdot,v),\derpar v
    C(\cdot,v)\right)\ud v\] in the class of all such $C\in \mathcal
  A$.
  \footnote{Note the difference between $\Len(C)$ and $\len(c)$, that was
    defined in eqn.~\eqref{eq:len}.
  }

  Any homotopy that achives the minimum of $\Len_F(C)$
  is called a \emph{geodesic}.
\end{Definition}
We call such distances $d_F(c_0,c_1)$, dropping the square brackets
for simplicity.%
\footnote{We are abusing notation:
these $d_F$ are not, properly speaking, distances in the space
$M$, since the distance between $c$ and a reparameterization 
$c\circ\phi$ is zero.}

\smallskip

We provide an interesting example of the above ideas.

\subsection{$L^\infinity$-type Finsler metric and Fr\'echet distance}
We digress from the main theme of the paper to prove a result 
that will be used in the following.

For any fixed immersed curve $c$ and $\theta\in S^1$,
we define for convenience $\pi_N :\real^n\to \real^n$
to be  the projection on the space
$N(\theta)$ orthogonal to the tangent vector $D_s c(\theta)$,
\begin{equation}
  \label{eq:pi_N}
  \pi_{N(\theta)} w = w - \langle w , D_s c(\theta)\rangle D_s
  c(\theta) \qquad \forall w\in\real^n.
\end{equation}

Consider two immersed curves $c_0$ and $c_1$; the Fr\'echet distance $d_f$
(as found in \cite{Michor-Mumford}) is defined by
\begin{Definition}[Fr\'echet distance]
  \[ d_f(c_0,c_1) \defeq \inf_\phi\sup_u |c_1(\phi(u))-c_0(u)| \]
  where $u\in S^1$ and $\phi$ is chosen in the class of
  diffeomorphisms of $S^1$.
\end{Definition}
This is a well defined distance in the space $B_i$ (that is not, though,
complete w.r.t. this distance: its completion is the space of
Fr\'echet curves).

Another similar distance was defined
in \S2.2 in \cite{YM:metrics04} by a different approach,
using a Finsler metric:
\begin{Definition}[Finsler $L^\infinity$ metric]
  \label{def:L_infinity}
  If we wish to define a norm $F(c,\cdot)$ on $T_cM$ that is modeled
  on the norm of the Banach space $L^\infinity(S^1\to\real^n)$, we  define
  \[ F^\infinity(c,h) \defeq \|\pi_N h\|_{L^\infinity}=\sup_\theta|\pi_{N(\theta)}  h(\theta)|\]

  We define the
  distance   $d_\infinity(c_0,c_1)$  as in \ref{def: d F}.
\end{Definition}
Section  \S2.2.1 in \cite{YM:metrics04} discusses the relationship
between the distance $d_\infinity$ and the Hausdorff distance of compact sets;
we instead discuss the relationship between  
$d_f$ and $d_\infinity$: indeed we prove that $d_f=d_\infinity$.

\begin{Theorem}\label{thm:d f = d inf}
  $d_f=d_\infinity$.

  \begin{proof}
    Fix $c_0$ and $c_1$,  and define $\mathcal A$ as in \ref{def: d F}.

    We recall that $d_\infinity$ is also equal to the infimum of 
    \[d_\infinity(c_0,c_1)=\inf_{C\in\mathcal A}\int_0^1
    \sup_\theta\big|\derparr v C(\theta,v)\big|\ud v\]
    as well (the proof follows immediatly from prop.~3.10 in  \cite{YM:metrics04})
    
    Consider a homotopy $C=C(u,v)\in\mathcal A$ connecting the curve $c_0$ to
    a reparameterization $c_1\circ\phi$ of the curve $c_1$:
    \[ \sup_{u} | c_1(\phi(u))-c_0(u)|= \sup_{u} | C(u,1)-C(u,0)|=\]
    \[=\sup_{u} \left|\int_0^1 \derparr v { C}(u,v)\ud v\right|\le
    \int \sup_{u} \Big|\derparr v { C}(u,v)\Big|\ud v\] 
    so that $d_f\le d_\infinity$.

    On the other side, let
    \[ C^\phi(\theta,v)\defeq (1-v) c_0(\theta)+ v c_1(\phi(\theta))\] 
    be the linear interpolation: then
    \[\derparr v { C^\phi}(u,v) = c_1(\phi(u))-c_0(u)\]
    (that does not depend on $v$) so that
    \[\sup_{u} \left|\int_0^1 \derparr v { C^\phi}(u,v)\ud v\right|=
    \int \sup_{u} \Big|\derparr v { C^\phi}(u,v)\Big|\ud v\] 
    and then, for that particular homotopy $C^\phi$,
    \[\Len_\infinity (C^\phi) =  \sup_{u} | c_1(\phi(u))-c_0(u)| \]
    we compute the infimum of all possible choices of
    $\phi$ and get that
    \[d_\infinity(c_0,c_1) =     \inf_C \Len_\infinity (C)  \le 
    \inf_\phi \Len_\infinity (C^\phi) =
    \inf_\phi \sup_{u} | c_1(\phi(u))-c_0(u)| =d_f(c_0,c_1)\]
  \end{proof}
\end{Theorem}
The theorem holds as well if use 
orientation preserving diffeomorphism $\mbox{Diff}^+(S^1)$
both in the definition
of the Fr\'echet distance and in the definition
of $L^\infinity$.

\section{Sobolev-type $H^j$ metrics}

We start by generalizing the definition in \ref{defn:inner}.
Fix $\lambda>0$.
Suppose that $h\in L^2$, then we can express it in Fourier series:
\begin{equation}
  \label{eq:hseries}
  h(s) = \sum_{l \in \Z} \hh(l) \exp{\left(\frac{2\pi i}{L}ls\right)} \\
\end{equation}
where $\hh  \in \ell^2(\Z\to\mathbb C)$; and similarly for $k$.

For any $\alpha>0$, given the Fourier coefficients $\hh,\hk : \Z \to
\mathbb{C}$ of $h,k$, we define the fractional Sobolev inner product
\begin{equation}
  \ip{h}{k}{H^\alpha_0}\defeq\sum_{l\in \Z}(2\pi l )^{2\alpha} \hh(l) \cdot \overline{\hk(l)}
\label{eq:frac_seminorm}    
\end{equation}
that is independent of curve scaling;
 then we can define
\[\ip{h}{k}{{H}^\alpha} \defeq \mean{h \cdot k}+\lambda 
\ip{h}{k}{H^\alpha_0}\]
\begin{equation}
  \ip{h}{k}{\tilde{H}^\alpha} \defeq \mean{h} \cdot \mean{k}+\lambda 
  \ip{h}{k}{H^\alpha_0}
  \label{eq:decompose_H_j_norm}
\end{equation}

When $\alpha=j$ integer,
these definition coincide with the one in \ref{defn:inner}.
So, for any $\alpha>0$, we represent the Sobolev--type metrics by
\begin{align}
  \label{eq:freqH1ip}
  \ip{h}{k}{H^\alpha} &= 
  \sum_{l\in \Z} (1+\lambda(2\pi l)^{2\alpha}) \hh(l) \cdot 
    \overline{\hk(l)} \\
    \label{eq:freqH1tip}
    \ip{h}{k}{\tilde{H}^\alpha} &= 
    \hh(0)\cdot \overline{\hk(0)} + \sum_{l\in \Z} \lambda
    (2\pi l)^{2\alpha} \hh(l) \cdot \overline{\hk(l)}.
  \end{align}

  \begin{Remark}
    Unfortunately for $j$ that is not an integer, the inner
    products (therefore, norms) are not local, that is, they cannot be
    written as integrals of derivatives of the curves.
    An interesting representation
    is by kernel convolution: given $r\in    \R^+$,
    we can represent them, for $j$ integer $j>r+1/4$, as
    \[
    \ip{h}{k}{\tilde{H}^r} = \int_{c} \int_{c}
    D^jh(s) K(s-\tilde s) D^jk(\tilde    s) \ud s \ud \tilde s
    \]
    that is,
    $
    \ip{h}{k}{\tilde{H}^r} = \ip{D^jh}{K\ast D^jk}{H^0}
    $,
    for a specific kernel $K$, ( $\ast$ denotes convolution in $S^1$
    w.r.t. arc parameter).
\end{Remark}

\begin{Remark}%
  \label{rem:transl_deform_decomp}
  The norm $\|h\|_{\tilde{H}^j}$ has an interesting interpretation
  in connection with applications in Computer Vision.
 
  Consider a deformation $h\in T_cM$ and write it as $h=\mean{h}+\tilde h$:
  this decomposes 
  \begin{equation}
    T_cM = \real^n \oplus D_cM\label{eq:transl_deform_decomp}    
  \end{equation}
  with
  \[D_cM\defeq \Big\{h:S^1\to\real^n ~|~ \mean{h}=0 \Big\}\]

  If we assign to $\real^n$ its usual euclidean norm, and to $D_cM$
  the scale-invariant $H^\alpha_0$ norm,  then we are naturally lead 
  to decompose as in eqn.~\eqref{eq:decompose_H_j_norm}, that is
  \begin{equation}
    \|h\|_{\tilde{H}^\alpha} ^2 = |\mean{h}|^2_{\real^n }+\lambda
    \|\tilde h\|_{{H}^\alpha_0}^2\label{eq:decompose_H_j}
  \end{equation}
  This means that the two spaces $ \real^n$ and $ D_cM$
  are orthogonal w.r.t. ${\tilde{H}^\alpha}$.

  In the above, $\real^n$ is akin to be  \emph{the space of translations} and
  $D_cM$ \emph{the space of non-translating deformations}.  That labeling 
  is not rigorous, though! since the subspace of $T_cM$ that does not move
  the center of mass $\mean c$ is not $D_cM$, but rather
  \[ \Big\{ h :  \int_{S^1}   h  + \big( c -\mean{c} \big) 
      \langle D_s h \cdot T \rangle \ud s = 0 \Big\} \]
  according to eqn.~\eqref{eq:der avg c}.

  The decomposition \eqref{eq:decompose_H_j} is at 
  the base of a two step alternating algorithm for minimization
  of $\tilde H^1$ Sobolev Active Contours (see \S3.4 in
  \cite{SYM05:sobol_activ_contour}),
  where the tracking of contours is done by, alternatively, minimizing
  an energy on curves translations, and then on curves deformations in
  $D_cM$ with the metric $\|\tilde h\|_{{H}^1_0}$. The resulting
  algorithm is independent of the choice of $\lambda$.

  Note that $\sqrt{\ip{h}{h}{H^\alpha_0}}$ is a norm on $D_cM$
  (by \eqref{eq:Poincare L2}), and
  it is a seminorm and not a norm on $T_cM$.
\end{Remark}

We define norms as
  \[F_{ H^j}(c,h)=
  \|h\|_{H^j} = \sqrt{\langle h,h\rangle_ {H^j}}
\quad,\quad
 F_{\tilde H^j}(c,h)=\|h\|_{\tilde{H}^j} = 
\sqrt{\langle h,h\rangle _{\tilde H^j}}
  \]
and consequently we define distances $d_{ H^j}$ and $d_{\tilde H^j}$
as explained in \ref{def: d F}.

\subsection{Stuff}

We improve a result from \cite{SYM05:sobol_activ_contour}
\footnote{and we provide a better version that unfortunately was
  prepared  too late
  for the printed version of \cite{SYM05:sobol_activ_contour}}:
we show that the norms associated with the inner products $H^j$ and
$\tilde{H}^j$  are equivalent.   We first prove 
\begin{Lemma}[Poincar\'e  inequalities]\def\l{\quad\Rightarrow\quad}
  Pick $h:[0,L]\to\real^n$,  weakly differentiable,
  with $h(0)=h(L)$ (so $h$  is periodically extensible) so that
  \[    h(u)-h(0)=\int_0^u h'(s) \ud s = -\int_u^L h'(s) \ud s    \]
  then derive these equations
  \begin{eqnarray*}
    h(u)-h(0)&=&
    \frac 1 2\left(\int_0^u h'(s) \ud s -\int_u^L h'(s) \ud s\right )  \l \\
    \l \mean h -h(0) &=&  \frac 1 {2L}
    \int_0^L \left(\int_0^u h'(s) \ud s -\int_u^L h'(s) \ud s\right) \ud u \l \\
    \l |\mean h -h(0)| &\le &  \frac 1 {2L}
    \int_0^L \left(\int_0^u |h'(s)|\ud s +\int_u^L |h'(s)|\ud
      s\right)\ud u = \\
    &=&\frac 1 {2L}
    \int_0^L \left(\int_0^L |h'(s)|\ud s \right)\ud u =
    \frac 1 2 \int_0^L |h'(s)| \ud s    
  \end{eqnarray*}
  so that (by  extending $h$ and replacing $0$ with an arbitrary point)
  \begin{equation}
    \sup_u |h(u)-\mean h| \le \frac 1 2 \int_0^L |h'(s)| \ud s 
    \label{eq:base_Poincare}    
  \end{equation} 
  \def\c{{\mathbb 1}}%
  the constant $1/2$ is optimal and is approximated by a family of $h$
  such that  $h'(s)=a(\c_{[0,\e)}(s)-\c_{[\e,2\e)}(s))$ when $\e\to 0$ (for
  a fixed $a\in\real^n$).%
  \footnote{$\c_A(x)$ is the characteristic function, taking value
    $1$ for $x\in A$,   $0$ for $x\not\in A$.  }
  
  By using H\"older inequality we can then derive many useful 
  Poincar\'e  inequalities of the
  form $\|h-\mean{h}\|_p\le c_{p,q,j} \|h'\|_q$.
  By Fourier transform we can also prove for $p=q=2$ that
  \begin{equation}
    {\int_0^L |h(s)-\mean h|^2 \ud s } \le
    \frac{L^{2j}}{(2\pi)^{2j}} {\int_0^L |h^{(j)}(s)|^2 \ud s}
    \label{eq:Poincare L2}    
  \end{equation}
  where the constant $c_{2,2,j}=(L/2\pi)^{2j}$
  is optimal and is achieved by $h(s)=a \sin(2\pi s /L)$ 
  (with $a\in\real^n$).\footnote{This $h$ is not the only solution;
    for $n=2$ we also have $h(s)=(\cos(2\pi s/L),\sin(2\pi s/L))$.  }
\end{Lemma}
\begin{Proposition}\label{prop:H1_equiv}
 Fix a smooth immersed curve $c:S^1\to\real^n$,
  let $L=\len(c)$.
 By H\"older's inequality, we have that
  $ |\mean{h}|^2 \le  \frac 1 L \int_0^L |h(s)|^2 \ud s $
  so that $\|h\|_{\tilde{H}^j}\le \|h\|_{H^j}$. On the other hand,  
  \begin{equation}
    \frac 1 L \int_0^L  |h(s)-\mean h |^2 \ud s =
    \frac 1 L \int_0^L  |h(s)|^2 \ud s - |\mean h |^2   \label{eq:variance}
  \end{equation}
  so that (by the  Poincar\'e  inequality \eqref{eq:Poincare L2}),
  \begin{eqnarray*}
    \|h\|_{H^j}^2 &=&
    \int_0^L \frac 1 L |h(s)|^2  +\lambda L^{2j-1} |h^{(j)}(s)|^2 \ud s \\
    &=& \frac 1 L \int_0^L  |h(s)-\mean h |^2 \ud s    + 
    \int_0^L \lambda L^{2j-1} |h^{(j)}(s)|^2  \ud s + |\mean h |^2 \\
    &\le&
    |\mean h |^2  +  L^{2j-1}\Big(\frac 1{(2\pi)^{2j}}+\lambda\Big)
    {\int_0^L |h^{(j)}(s)|^2 \ud s} \le
    \frac{1+{(2\pi)^{2j}}\lambda }{(2\pi)^{2j}\lambda} \|h\|_{\tilde{H}^j}^2
  \end{eqnarray*}
\end{Proposition}
Consequently,
\[ d_{ \tilde H^j}\le  d_{  H^j}\le 
\sqrt{\frac{1+(2\pi)^{2j}\lambda }{(2\pi)^{2j}\lambda} }
d_{ \tilde H^j}\]

More in general
\begin{Proposition}\label{prop:Hj_equivalent}
  For $i=0,\ldots,j$,
  choose $\overbar a_0\ge 0$ and $a_i\ge 0$ with $a_0+\overbar a_0>0$ and
  $a_j>0$.
  Define a $H^j$-type Riemannian norm~%
  \footnote{the scalar product can  be easily inferred}
  \begin{equation}
    \|h\|_{(a),j}^2 \defeq\overbar a_0|\mean{h}|^2+
    \sum_{i=0}^j a_i L^{2i-1} \int_0^L |h^{(i)}(s)|^2ds  
    \label{eq:Hj_general}    
  \end{equation}
  then all such norms are equivalent.

  \def\j{r}
  Moreover, choose $ \j$ with $1\le \j\le j$, and choose
  $\overbar b_0\ge 0$,
  $b_i\ge 0$ with $\overbar b_0+ b_0>0,b_{\j}>0$: then
  the norm $\|h\|_{(a),j}$ is stronger than the norm  $\|h\|_{(b),\j}$.

  \begin{proof}
    The proof is just an application of \eqref{eq:variance} and of
    \eqref{eq:Poincare L2} (repeatedly); note also that for $1\le i<j$
    equation~\eqref{eq:Poincare L2} becomes
    \begin{equation}
      {\int_0^L |h^{(i)}(s)|^2 \ud s } \le
      \frac{L^{2j-2i}}{(2\pi)^{2j-2i}} {\int_0^L |h^{(j)}(s)|^2 \ud s}
      \label{eq:Poincare L2 j i}    
    \end{equation}
    since $\mean{h^{(i)}}=0$.
  \end{proof}
\end{Proposition}

So our definitions of  $\|\cdot\|_{H^j}$ and  $\|\cdot\|_{\tilde H^j}$ are
in a sense the simpler choices of a Sobolev type norm that are scale invariant;
in particular,
\begin{Remark}\label{MM equiv}
  the $H^j$ type metric 
  \[\|h\|_M^2\defeq \int \sum_{i=0}^j |h^{(i)}(s)|^2 \ud s\]
  studied  in \cite{MM:hamil_rieman} is equivalent to our choices,
  \[b_1\|\cdot\|_{\tilde H^j }\le  \|\cdot\|_M \le b_2\|\cdot\|_{\tilde H^j} \]
  but the constants $b_1,b_2$ depend on the length of the curve.
\end{Remark}

\medskip

Following this proposition, we will prove some properties of the $H^1$ 
metric, and we will know that they can be extended to  $\tilde H^1$ and to
more general $H^j$-type metrics defined as  in \eqref{eq:Hj_general}.

\bigskip

We prove this fundamental inequality~\eqref{eq:H1_funda}:
\begin{Proposition} Fix a smooth immersed curve $c:S^1\to\real^n$,
  let $L=\len(c)$.

  We rewrite for convenience
  \[ \|h\|_{H^1}^2\ge \lambda L^2 \ip{h'}{h'}{H^0}= \lambda L \int_0^L
  |h'(s)|^2 \ud s = \lambda \int |\dot c(u)|\ud u \, \int
  |h'(u)|^2 |\dot c(u)| \ud u  \]
  where $h'=D_s h$;
  then by Cauchy-Schwartz
  \[ \int |\dot c(u)|\ud u \, \int |h'(u)|^2 |\dot c(u)| \ud u \ge
  \left(\int |h'(u)| |\dot c(u)| \ud u\right)^2 \]
  Suppose now that
  $C(u,v)$ is a smooth  homotopy of immersed curves $C(\cdot,v)$:
  then set $h(u,v)=\derpar v C(u,v)$ so that
  $D_s h=D_s \derpar v C=\derpar {uv} C /|\derpar u C|$. Summarizing the above 
  \begin{equation}
    \label{eq:H1_funda}
    \| \derpar v C(\cdot,v)\|_{H^1}\ge 
    \sqrt\lambda \int  |\derpar {uv}C(u,v)| \ud u
  \end{equation}

  As argued in \ref{prop:Hj_equivalent}, the above result extends to 
  all $H^j$-type norms \eqref{eq:Hj_general} .
\end{Proposition}

We related the $H^1$-type metric to the $L^\infinity$ type metrics
\begin{Proposition}\label{prop:H_and_Frechet}
  The $\tilde H^1$ metric is stronger than the $L^\infinity$ metric defined in 
  \ref{def:L_infinity}.

  As a consequence, by theorem \ref{prop:H1_equiv} and \ref{thm:d f = d inf},
  the $H^j$ and $\tilde H^j$ distances 
  are lower bounded by the  Fr\'echet distance
  (with appropriate constants depending on $\lambda$).
 
  \begin{proof}
    Indeed, by \eqref{eq:base_Poincare}    there follows
    \begin{eqnarray*}
      \sup_\theta|\pi_{N(\theta)}  h(\theta)| &\le&
    \sup_\theta| h(\theta)| \le
    |\mean{h}| +\frac 1 2 \int |h'|ds\le \\  &\le&
    |\mean{h}| +\frac{\sqrt L} 2 \sqrt{\int |h'|^2 ds}\le 
    \sqrt 2\sqrt{|\mean{h}|^2 +\frac{L} 4 \int |h'|^2 ds}
  \end{eqnarray*}
  ($\pi_N$ was defined in eqn.~\eqref{eq:pi_N}). For example,
  choosing $\lambda=1/4$,  
  \[ F_\infinity (c,h) \le \sqrt 2 \|h\|_{\tilde H^1}\]
  \end{proof}
\end{Proposition}

We also establish relationship between the length $\len(c)$ of a curve
and the Sobolev metrics:
\begin{Proposition}\label{prop:H1_and_len}
  Suppose again that
  $C(u,v)$ is a smooth homothopy of immersed curves,
  let $L(v)\defeq\len(C(\cdot,v))$ be the length at time $v$; then  
  \[ \derpar v L = 
  \int \langle \derpar {uv} C,\frac{\derpar u  C}{|\derpar u  C|}\rangle \ud u
  \le   \int |\derpar {uv} C| \ud u \le
 \frac 1{\sqrt\lambda}   \| C_v(\cdot,v)\|_{H^1}\]
  by \eqref{eq:H1_funda}.
\end{Proposition}

We have many interesting consequences:
\begin{itemize}
\item 
  \begin{equation}
    |L(1)-L(0)| \le \frac 1{\sqrt\lambda} \Len(C)\label{eq:L_Lip}      
  \end{equation}
  where the length $\Len(C)$ of the homotopy/path $C$ is computed
  using either $H^1$ or $\tilde H^1$ (or using any metric as in
  \eqref{eq:Hj_general} above, but in this case the constant in
  \eqref{eq:L_Lip} would change).
\item 
  Define the length functional $c\mapsto \text{len}(c)$
  on our space of curves;  embed the space of curves with a  $H^1$ metric;
  then  the length functional is Lipschitz.
  
\item 
  The ``zero curves'' are the constant curves (that have zero
  length); these are points in the space of curves where the space
  of curves is, in a sense, singular; by the above, the ``zero
  curves'' are a closed set in the $H^1$ space of curves, and an immersed
  curve $c$ is distant at least $\len(c){\sqrt\lambda}$ 
  from the ``zero curves''.
\end{itemize}

But the most  interesting consequence is that
\begin{Theorem}[Completion of $B_1$ w.r.t. $H^1$]
  \label{prop:H1_completion}
  let  $d_{ H^1}$  be the distance induced by $H^1$; 
  the metric completion of the space of curves is contained in the space of all
  rectifiable curves.
\end{Theorem}
\begin{proof}
  This statement is a bit fuzzy: indeed $d_{ H^1}$ is
  not a distance on $M$, whereas in $B$ objects are not functions, but
  classes of functions.  So it must be intended ``up to
  reparametrization of curves'', as follows.%
  \footnote{The concept is clarified by introducing the concept
    of \emph{horizontality} in $M$, that we must unfortunately skip for sake
    of brevity  }

  Let $(c_n)_{n\in\natural}$ be a Cauchy sequence.  Since $d_{ H^1}$
  does not depend on parametrization, we assume that all $c_n$ are
  parametrized by arc parameter, that is, $|\derpar \theta c_n|=l_n$
  constant in $\theta$.  By proposition \ref{prop:H_and_Frechet} , all
  curves are contained in a bounded region; since $\len(c_n)=2\pi l_n$
  by proposition \ref{prop:H1_and_len} above,  the sequence $l_n$
  is bounded. So the (reparametrized) family $(c_n)$ is equibounded
  and equilipschitz: by Ascoli-Arzel\`a theorem, up to a subsequence,
  we obtain that $c_n$ converges uniformly to a
  Lipschitz curve $c$, and  $|\derpar \theta c|\le \lim_n l_n$.
\end{proof}
We also prove that 
\begin{Theorem}\label{thm:H1_approx} 
Any rectifiable planar curve is approximable by smooth curves
according to  the distance induced by $H^1$.
\end{Theorem}
\begin{proof}
  Let $c$ be a rectifiable curve, and assume that it is    non-constant.

    As a first step, we assume that $c$ is not \emph{flat}, that is,
    the image of $c$ is not contained     in a line in the plane.
  We sketch how we can approximate $c$    by smooth curves.
    The precise arguments are in section 2.1.4 in
    \cite{YM:metrics04}; see in particular the proofs
    of \ref{03-prop: Sriva M} and \ref{03-prop: no rotation Srivastava}
    in the appendix.  Since
    the metric is independent of rescaling, we rescale $c$, and assume
    that $|\derpar \theta c|=1$.  

    We identify $S^1$ with $[0,2\pi)$. Let in the following $L^2= L^2([0,2\pi]$.
    We define the \emph{measurable angle function} to be a
    function $\tau:[0,2\pi)\to [0,2\pi) $ such that $\derpar \theta
    c(\theta) = (\cos \tau(\theta),\sin \tau(\theta))$.  
    We define
    \[ S=\left\{ \tau\in L^2([0,2\pi]) ~|~ \phi(\tau)=(0,0)
    \right\}\] where $\phi:L^2\to\real^2$ is defined by
    \[  \phi_1(\tau) =    \int_0^{2\pi}\!\cos\tau(s) ds ~~,~~ \phi_2(\tau) =
    \int_0^{2\pi}\!\sin\tau(s) ds \] 
    (this is similar to what was done
    in Srivastava et al. works on \emph{``Shape Representation using
      Direction Functions''}, see \cite{Srivastava:AnalPlan}).
    
    As proved in \ref{03-prop: Sriva M} in \cite{YM:metrics04},
    $S$ is a manifold near $\tau$. 
    As shown in the proof of \ref{03-prop: no rotation Srivastava}
    in \cite{YM:metrics04}, there exists a function $\pi:V\to S$ defined 
    in a neighbourhood $V\subset L^2$ of $\tau$ such that,
    if $f(s)\in L^2$ is smooth, then $\pi(f)(s)$ is smooth.
    Let $f_n$ be a smooth approximation of $\tau$, with $f_n\to \tau$
    in $L^2$; then $g_n\defeq \pi(f_n)\to \tau$.  Let then
    \[ G_n(\theta,t) \defeq \pi(t \tau + (1-t) g_n )(\theta)  \]
    be a the projection on $S$ of the linear path connecting $\tau$ to
    $g_n$. Since $S$ is smooth in $V$, then
    the $L^2$ distance $\|\tau - g_n \|$ is equivalent to the
    geodesic induced distance; in particular, 
    \[\lim_n \E_S G_n = 0 \]
    where
    \[\E_S G \defeq \int_0^1 \|\derpar t G(\cdot ,t)\|_{L^2}^2 \ud t \]
    is the action of the path $G$ in $S\subset L^2$.

    The above $G_n$ can be associated to an homotopy by definining
    \[ C_n(s,t) \defeq c(0) +
    \int_0^s (\cos (G_n(\theta ,t)),\sin(G_n(\theta ,t)) )\ud \theta\]
    note that $C_n(s,0)=c(s)$ and $C_n(s,1)$ is a smooth closed curve.

    We now compute the $\tilde H^1$ action of $C_n$,
    \[\E_{\tilde H^1}(C_n) \defeq \int_0^1 \| \derpar t C_n\| _{H^1}^2\ud t =
    \int_0^1 %
      \int_0^{2\pi } | \derpar t C_n|^2 +
      | D_s \derpar t C_n|^2\ud s \ud t
    \]
    Since any $C_n(\cdot,t)$ is by arc parameter, then 
    $  D_s \derpar t C_n=  \derpar{s t} C_n$ so
    \[ D_s \derpar t C_n =    N(s) \derpar t G_n(s ,t)    \]
    where 
    \[N(s)\defeq (-\sin (G_n(s ,t)),\cos(G_n(s ,t)) ) \]
    is the normal to the curve; so the second term 
    in the action $\E_{ H^1}(C_n)$ is exactly equal to 
    $\E_{S}(G_n)$, that is,
    \[\E_{ H^1}(C_n) = \int_0^1 \| \derpar t C_n\| _{H^1}^2\ud t =
    \int_0^1 
      \int_0^{2\pi } | \derpar t C_n|^2 \ud s \ud t + E_{S}(G_n)    \]
    We can also prove that $\int_0^1
      \int_0^{2\pi } | \derpar t C_n|^2 \ud s \ud t\to 0$, so 
      $\E_{ H^1}(C_n)\to 0$, and then
    \[\lim_n \Len_{ H^1} (C_n)=0    \]

    \medskip 

    As a second step, to conclude,
    we assume that $c$ is  \emph{flat}, that is,
    the image of $c$ is  contained     in a line in the plane;
    then, up to translation and rotation,
    \[ c(\theta) = (c_1(\theta),0)\]
    since $c$ is by arc parameter, $\dot c_1 = \pm 1$.
    Let then $f:[0,2\pi]$ be smooth and with support in $[1,3]$ and
    $f(2)=1$; let    moreover
    \[ C(\theta,t) \defeq ( c_1(\theta) ,t f(\theta)) \]
    so
    \[ |\derpar \theta C| = \sqrt{1 + (f'(\theta))^2} \ge 1    \]
    and then we can easily prove that 
    \[\Len_{ H^1} (C)<\infinity    \]
    moreover, any curve $C(\cdot,t)$ for $t>0$ is not flat,
    so it can be approximated by smooth curves
\end{proof}

\subsection{The completion of  $M$ according to $H^2$ distance}

Let $d(c_0,c_1)$ be the geometric distance induced by  $H^2$ on $M$
(as defined in \ref{def: d F}).
Let $E(c)\defeq \int|D_s^2 c|^2\ud s$ be defined on non-constant
smooth curves.
We prove that
\begin{Theorem}\label{thm:completionH2}
  $E$ is locally Lipschitz in $M$ w.r.t. $d$, and the local Lipschitz
  constant depends on the length of $c$. 

  As a corollary, all
  non-constant curves in the completion of
  $C^\infinity(S^1\to\real^n)$
  according to  the metric $H^2$ admit curvature as a measurable function,
  and the curvature satisfies $E(c)<\infinity$.
  
  Viceversa, any non-constant curve admitting curvature in a weak
  sense and satisfying $E(c)<\infinity$ is approximable by smooth curves.
\end{Theorem}

The rest of this section is devoted to proving the above three statements.

Fix a curve $c_0$; let $L_0\defeq \len c_0$ be its length.

By eqn.~\eqref{eq:L_Lip} and prop.~\ref{prop:Hj_equivalent},
we know that the ``length function'' 
$c\mapsto\len(c)$ is Lipschitz in $M$ w.r.t the distance $d$,
that is,
\[ |\len c_0 - \len c_1|\le a_1 d(c_0,c_1) \]
where $a_1$ is a positive constant (dependent on $\lambda$).

Choose any $c_1$ with $d(c_0,c_1)< L_0/(4a_1)$.

Let $C(\theta,t)$ be a time varying smooth homotopy
connecting $c_0$ to (a reparametrization of) $c_1$; choose it so that
$\Len C < 2 d(c_0,c_1) $; then $\Len C < L_0/(2a_1)$.

Let $L(t)\defeq \len C(\cdot,t)$ be the length of the curve at time $t$.
Since at all times $t\in[0,1]$, $d(c_0,C(\cdot,t))< L_0/(2a_1)$,
then  $|L(t)- L_0|< a_1 L_0/(2a_1)=L_0/2$; in particular,
\[ L_0/2 < L(t) < L_0 3/2 ~ .\]
By using this last inequality,
we are allowed to discard $L(t)$ in most of the following estimates.

We call 
$ \| f \| \defeq \sqrt{\int |f(s)|^2\ud s}$
and
\[ N(t) \defeq
 \| D_s^2 \derpar t C(\cdot,t) \|=\sqrt{\int |D_s^2 \derpar t C|^2 \ud s } \]
for convenience; using this notation, we recall that 
\[  \| \derpar t C \|_{H^2} =  \sqrt{  \lambda L(t)^3 N(t)^2
+ \frac 1 {L(t)}  \|\derpar t C \|^2 } ~~;\]
so $\| \derpar t C \|_{H^2}\ge \sqrt \lambda L^{3/2} N(t)$.

Up to reparametrization in the $t$ parameter, we can suppose that the path
$t\mapsto C(\cdot,t)$ in $M$ is by (approximate) arc parameter,
 that is $\| \derpar t C \|_{H^2}$ is (almost) constant in $t$;
so we assume, with no loss of generality, that
$ \| \derpar t C \|_{H^2} \le 2 d(c_0,c_1)$ for all $t\in[0,1]$, and then
$  N(t)  \le a_2 d(c_0,c_1)$ where 
$a_2 = 2 / \sqrt{ (L_0/2) ^{3}  \lambda}$.

We want to prove that
\[ E(c_1) - E(c_0) \le a_5 d(c_0,c_1)  \]
where the constant $a_5$ will depend on $L_0$ and $\lambda$.

By direct computation
\[ \derpar t E(C(\cdot,t)) =
 \int|D_s^2 C|^2 \langle D_s \derpar t  C ,D_s C\rangle\ud s +
2\int \langle D_s^2 C, \derpar t  D_s^2 C\rangle \ud s\]

we deal with the two addenda in this way:
\begin{enumerate}
\item by Poincar\'e inequality \eqref{eq:base_Poincare} we deduce
  \[\sup_\theta
  |D_s \derpar t C| \le \frac 1 2 \int |D_s^2 \derpar t C| \ud s \le
  \sqrt{L(t)}\sqrt{\int |D_s^2 \derpar t C|^2 \ud s} =
  \sqrt{L(t)} {N(t)}\]
  since $\mean{D_s \derpar t C}=0$.

  So we estimate the first term as
  \[ \int|D_s^2 C|^2 \langle D_s \derpar t  C ,D_s C\rangle\ud s
  \le E(C) \sqrt{L(t)} {N(t)} ~. \]

\item
 The commutator of $D_s$ and $\derpar t$ is 
 ${\langle D_s \derpar t c, D_s c  \rangle} D_s$: indeed
  \begin{eqnarray*}
   \derpar t  D_s &=&
  \frac 1 {|\derpar \theta c|} \derpar \theta  \derpar t  +  
  (\derpar t \frac 1 {|\derpar \theta c|} ) \derpar \theta  =
   D_s  \derpar t - 
 \frac {\langle\derpar t \derpar\theta c,\derpar\theta c\rangle} {|\derpar \theta c|^3}\derpar \theta
  = \\ & = & D_s  \derpar t  - 
  {\langle D_s \derpar t c, D_s c  \rangle} D_s
\end{eqnarray*}

so 
 \begin{eqnarray*}  \derpar t  D_s^2 C &=&  D_s \derpar t D_s C - 
{\langle D_s \derpar t C, D_s C  \rangle}D_s^2 C
= \\ &=&
 D_s^2 \derpar t  C -
D_s({\langle D_s \derpar t C, D_s C  \rangle}D_s C) -
{\langle D_s \derpar t C, D_s C  \rangle}D_s^2 C = \\ &=&
 D_s^2 \derpar t  C -
({\langle D_s^2 \derpar t C, D_s C  \rangle}D_s C) -
({\langle D_s \derpar t C, D_s^2 C  \rangle}D_s C) - \\ & & -
2({\langle D_s \derpar t C, D_s C  \rangle}D_s^2 C)
\end{eqnarray*}
so (since $|D_s C|= 1$)
\[ \| \derpar t  D_s^2 C \| \le
 2 \| D_s^2 \derpar t  C\| + 3 \| D_s^2 C\| \sup|D_s\derpar t C|\]
that yields an estimate of the second term
\[\int \langle D_s^2 C, \derpar t  D_s^2 C\rangle \ud s\le 
\sqrt{E(C)} \Big( 2 N(t) + 2 \sqrt{E(C)}\sqrt{L(t)} {N(t)}  \Big)\]
by using Cauchy-Schwartz.
\end{enumerate}

Summing up
\[ |\derpar t E(C(\cdot,t))| \le 
2 \sqrt{E(C)}  N(t) + 3 {E(C)}\sqrt{L(t)} {N(t)} \]
or, since  $\sqrt{x}\le 1+x$,
\[ |\derpar t E(C(\cdot,t))| \le 
2 N(t)+2 {E(C)}  N(t)+ 3 {E(C)}\sqrt{L(t)} {N(t)} \]

We recall that $N(t)  \le a_2 d(c_0,c_1)$, $L(t)\le L_0 3 /2$,
so we rewrite the above as
\[ |\derpar t E(C(\cdot,t))| \le 
2  a_2 d(c_0,c_1) +   2 E(C)   a_2 d(c_0,c_1)   + 3 E(C)  a_4 a_2 d(c_0,c_1) \]
with $a_4=\sqrt{L_0 3 / 2}$.
Apply Gronwall's Lemma to obtain
\[ E(c_1) \le \Big( E(c_0)  + 2 a_2 d(c_0,c_1) \Big) 
\exp \Big(  (2+3 a_4) a_2 d(c_0,c_1)  \Big)  ~~. \]

Let
\[ g(y) \defeq  \Big( E(c_0)  + 2 a_2 y \Big) 
\exp \Big(  (2+3 a_4) a_2 y  \Big) \]
then $E(c_1)\le g(d(c_0,c_1))$; since $g$ is convex, and $g(0)= E(c_0)$,
then there exists a $a_5>0$ such that  
$g(y) \le  E(c_0)  +  a_5 y$ when $0\le y\le L_0/(4a_1)$;
since we assumed that  $d(c_0,c_1)< L_0/(4a_1)$, then 
\[ E(c_1) \le  E(c_0)  +  a_5 d(c_0,c_1) ~. \]
Note that  $a_5$ is ultimately dependent on $L_0$ and $\lambda$.

This ends the proof of the first statement of \ref{thm:completionH2}.

\medskip

To prove the second statement, 
let $(c_n)_{n\in\natural}$ be a Cauchy sequence.  Since 
$d_{ H^1}\le a d_{ H^2}$, then as in the proof of
\ref{prop:H1_completion}, we assume that, up to reparametrization
and a choice of subsequence,  $c_n$ converges uniformly to a
Lipschitz curve $c$.  

Let $L_0=\len c$.
We have assumed in the statement that $c$ is non-constant; then $L_0>0$.

Again, the ``length function'' $c\mapsto\len(c)$ is Lipschitz, so we
know that the sequence $\len(c_n)$ is Cauchy in $\real$, so it
converges; moreover the ``length function'' $c\mapsto\len(c)$ is lower
semicontinuous w.r.t. uniform convergence, so $\lim_n\len(c_n)\ge
\len(c)>0$. So we assume, up to a subsequence, that $2L_0 \ge \len(c_n)
\ge L_0$

We proved above that, in a neighbourhood of $c$ of size
$L_0/(8a_1)$, the function $E(c)\defeq \int|D_s^2 c|^2\ud s$ is
Lipschitz; so we know that the sequence $E(c_n)$ is bounded, and then
(since curves are by arc parameter and $\len(c_n)\ge L_0$) the energy
$\int|\derpar \theta^2 c|^2\ud s$ is bounded: then
$\derpar \theta c_n $ are uniformly H\"older continuous, so 
by Ascoli-Arzel\`a compactness theorems,
up to a subsequence,  $\derpar\theta c_n(\theta)$
converges.

As a corollary we obtain that $\lim_n\len(c_n)=\len(c)$, that $c$ is
parametrized by arc parameter, and that $D_s c_n(\theta)$ converges to
$D_s c(\theta)$.

Since the functional $\int|\derpar\theta^2 c_n|^2\ud s$ is bounded in
$n$, then by a theorem in \cite{Brezis}, $c$ admits weak derivative 
$\derpar\theta^2 c$ and $\int|\derpar\theta^2 c|^2\ud s<\infinity$,
and equivalently,  $\int|D_s^2 c|^2\ud s<\infinity$.

\medskip

For the third statement, viceversa, let $c$ be a rectifiable curve,
and assume that it is non-constant, and $E(c)<\infinity$. 
Since the metric is independent of rescaling, we rescale
$c$, and assume that $|\derpar u c(u)|=1$.

We express in Fouries series
\begin{equation}
 c(u) = \sum_{n\in\integer} l_n \exp( i n   u)\label{eq:F_c}
\end{equation}
(by equating $S^1=\real/2\pi$),   then   we
decide that 
\begin{equation}
  C(u,t) \defeq  \sum_{n\in\integer}  l_n
   \exp( i n   u - f(n) t)\label{eq:F_C}
 \end{equation}
with $f(n)=f(-n)\ge 0$ and $\lim f(n) / \log(n) =\infinity$;
(for example, $f(n)=|n|$ or $f(n)=(\log(|n|+2))^2$): then
$C(\cdot,t)$ is smooth for any $t>0$

\def\C{{\tilde C}}
We want to prove that, for $t$ small, $C(\cdot,t)$ is near $c$
in the $H^2$ metric; to this end, let $\C$ be the linear interpolator
\begin{equation}
   \C(u,t,\tau) \defeq (1-\tau) c(u) + \tau C(u,t)
   = \sum_{n\in\integer}  l_n
   e^{ i n u} ( 1-\tau  +\tau e^{- f(n) t})
   \label{eq:F_tC}
 \end{equation}
we will prove that 
\begin{equation}
\int_0^1
\left(\int_{S^1} |\derpar \tau \C|^2 \ud s+
\lambda L^4 \int_{S^1} |D^2_s \derpar \tau \C|^2 \ud s
\right)
 \ud \tau < \delta(t)\label{eq:short}
\end{equation}
where $\lim_{t\to 0}\delta(t)=0$, and $L$ is the length of $\C(\cdot,t,\tau)$.

We need some preliminary results:
\begin{itemize}
\item we prove that 
  \begin{equation}
    \label{eq:d1}
    \int_{S^1}|\derpar{uu} c-\derpar{uu} \C|^2\ud u<\delta_1(t)
  \end{equation}
  where $\lim_{t\to 0}\delta_1(t)=0$, uniformly in $\tau\in[0,1]$; we write
    \[\int_{S^1}|\derpar{uu} c-\derpar{uu} \C|^2\ud u=
    2\pi \tau ^2 \sum_{n\in\integer}  |l_n|^2 |n|^4
    (1- e^{- f(n) t})^2    \]
    and since
    \[    2\pi \sum_{n\in\integer}  |l_n|^2 |n|^4=
    E(c)=\int_{S^1}|\derpar{uu} c|^2\ud s<\infinity       \]
    and $\lim_{t\to 0}(1- e^{- f(n) t})^2=0$,
    we can apply Lebesgue dominated convergence theorem.
  \item We prove that 
    \begin{equation}
      \label{eq:d2}
      |\derpar{u} c-\derpar{u} \C|<\delta_2(t)
    \end{equation}
    where $\lim_{t\to 0}\delta_2(t)=0$, uniformly in $u$ and 
    $\tau\in[0,1]$; indeed
    \begin{eqnarray*}
    |\derpar{u} c-\derpar{u} \C|\le 
     \tau  \sum_{n\in\integer}  |l_n| |n|    (1- e^{- f(n) t}) \le \\ \le
     \sqrt{\sum_{n\in\integer}  |l_n|^2 |n|^4}  
     \sqrt{\sum_{n\in\integer,n\neq 0} \frac 1{n^2}  (1- e^{- f(n) t})^2 }
   \end{eqnarray*}

     again we apply Lebesgue dominated convergence theorem.
  \item By the above we also obtain that for $t$ small, 
    \begin{equation}
    3/2\ge |\derpar{u} \C|\ge 1/2 \text{ uniformly in }\tau,u
    \label{eq:uni2}    
  \end{equation}
\item
  We can similarly prove that
  \begin{equation}
      \label{eq:d3}
      |c- \C|<\delta_3(t)
    \end{equation}
\end{itemize}

By direct computation
\[ D^2_s \derpar \tau \C = 
\frac {\derpar{u u  \tau}  \C}{|\derpar u  \C|^2}
+  \frac {\langle \derpar{u u } \C,\derpar{u } \C  \rangle
       \derpar{u  \tau}  \C}
{|\derpar u  \C|^4} \]

but then, for $t$ small, by \eqref{eq:uni2},
\[  |D^2_s \derpar \tau \C|  \le
4 |\derpar{u u  \tau}  \C| 
+ 24 |\derpar{u u } \C| | \derpar{u  \tau}  \C| \]

We use the fact that
  \[ \derpar{u u  \tau}  \C = \derpar{u u  }  C  - \derpar{u u  }c 
~~,~~ \derpar{u   \tau}  \C = \derpar{u  }  C  - \derpar{u   }c 
~~,~~ \derpar{\tau}  \C = C  - c ,
  \]
  so by eqn.~\eqref{eq:d1} and eqn.~\eqref{eq:d2}
  \[ \int \int  |D^2_s \derpar \tau \C|^2\ud s \ud \tau  
  \le  a_1(\delta_1(t)  +E(c)\delta_2(t) )\]
and by eqn.~\eqref{eq:d3}
$\int \int  |\derpar \tau \C|^2\ud s \ud \tau    \le  8 \delta_3(t)$.
Eventually we combine all above to  bound eqn.~\eqref{eq:short} by
setting $\delta(t) =a_2 \delta_3(t) + \lambda a_2 (\delta_1(t)
+E(c)\delta_2(t) ) $.

This concludes the proof.

\tableofcontents

\bibliographystyle{amsalpha}
\bibliography{mennbib.bib}

\end{document}